\documentclass[12pt,a4paper]{article}

\usepackage{nicefrac}
\usepackage[latin1]{inputenc}
\usepackage{amsmath}
\usepackage{amssymb}
\usepackage[english]{babel}
\usepackage[T1]{fontenc}

\newtheorem{lem}{Lemma}[section]

\newtheorem{thm}[lem]{Theorem}

\newtheorem{cor}[lem]{Corollary}

\newtheorem{rem}[lem]{Remark}

\numberwithin{equation}{section}

\usepackage{color}

\vfuzz \hfuzz
\begin{document}

\title{Artificial boundary conditions for  linearized
stationary incompressible viscous flow around rotating and translating body}
\date{}
\author{P. Deuring $^1$, S. Kra\v cmar $^{2,3}$, \v S. Ne\v casov\' a$^4$}
\maketitle
\date{}
\centerline{$^1$ Univ. Littoral C\^ote d'Opale, Laboratoire de math\'ematiques} \centerline{pures
et appliqu\'ees Joseph Liouville}
\centerline{e-mail: {\tt Paul.Deuring@lmpa.univ-littoral.fr}}

\vspace{2mm}
\centerline{$^2$ Department of Technical Mathematics, Czech Technical University}

\centerline{$^3$ Institute of Mathematics of the Academy of Sciences of the Czech Republic}
\centerline{e-mail: {\tt Stanislav.Kracmar@fs.cvut.cz}}

\vspace{2mm}
\centerline{$^4$ Institute of Mathematics of the Academy of Sciences of the Czech Republic}

\centerline{e-mail: {\tt matus@math.cas.cz}}

\bigskip

\date{}
\vskip0.25cm
\begin{abstract}
We consider the linearized and nonlinear stationary incompressible flow around rotating and translating body in the exterior domain $\mathbb R^3 \setminus \overline{\mathcal D}$, where $\mbox{$\mathcal D$} \subset \mathbb{R}^3 $ is open and bounded,
with Lipschitz boundary. We derive the pointwise estimates for the pressure in both cases. Moreover, we consider the linearized problem in a truncation domain
$\mathcal D_R:=B_R \backslash \overline{ \mathcal D}$
of the exterior domain
$\mathbb R^3 \setminus \overline{\mathcal D}$
under certain artificial boundary conditions on the truncating boundary $\partial B_R$,
and then compare this solution with the solution  in the exterior domain $\mathbb R^3 \setminus \overline{\mathcal D}$ to get the truncation error estimate.
\end{abstract}

\section{Introduction}

We consider the systems of equations
\begin{equation}\label{1.0}
\begin{array}{crl}
-\Delta u (z) + ( \tau e_1- \varrho e_1 \times z)\cdot \nabla u(z) + \varrho e_1 \times u (z)
\ \ \ \ \ \ \ \ \ \ \ \ \ \ \ \ \ \ \ \ \ \\ \qquad\qquad\qquad\qquad
+\tau(u(z)\cdot \nabla) u(z)+\nabla \pi(z)  = F(z)\\
{\rm div}\,u(z)=0 \ \text{for} \ z\in \mathbb R^3 \setminus \overline{\mathcal D}
\end{array}
\end{equation}

\begin{equation}\label{eq_1.1}
\begin{array}{crl}
-\Delta u (z) + ( \tau e_1- \varrho e_1 \times z)\cdot \nabla u(z) + \varrho e_1 \times u (z)
+\nabla \pi(z)  = F(z)\\
{\rm div}\,u(z)=0 \ \text{for} \ z\in \mathbb R^3 \setminus \overline{\mathcal D}\\

\end{array}
\end{equation}where $\mbox{$\mathcal D$} \subset \mathbb{R}^3 $ is open and bounded,
with Lipschitz boundary.
Problems (\ref{1.0}) and (\ref{eq_1.1}) together with some boundary conditions on $ \partial \mathcal D $ constitute  mathematical models (linear  and non-linear, respectively) describing stationary flow of a viscous incompressible fluid around a rigid body which moves at a constant velocity and rotates at a constant angular velocity, where we consider that the rotation is parallel to the velocity at the infinity. { For details concerning of deriving the model, see \cite{Fa2, G2}. The description and the analysis in the case when the rotation is not parallel to the velocity at infinity  can be find in the following works, see \cite{FGTN, GH}}.

 The aim of this paper is two folds:

First, we would like to derive  the pointwise estimates for the pressure in the linear
 and
also in the non-linear cases in
order to complete the pointwise estimates for the velocity  and its gradient from \cite{ DKN5,3} by the pointwise estimates of the pressure in order to get complete decay  information  of all parts $u,\,\pi  $ of  solutions to systems (\ref{1.0}), (\ref{eq_1.1}).
{ Let us mention  that the decay of pressure was also investigated  in the work of Galdi, Kyed \cite{GK1} and in case of pure rotation  see \cite{FGK}}.

 Second, to solve the linear system (\ref{eq_1.1}) in a truncation
$\mathcal D_R:=B_R \backslash \overline{ \mathcal D}$
of the exterior domain
$\mathbb R^3 \setminus \overline{\mathcal D}$
under certain artificial boundary conditions on the truncating boundary $\partial B_R$,
and then compare this solution with the solution of (\ref{eq_1.1}) in the exterior domain, i.e. to get some sort of error estimates of the method of an artificial boundary condition. For this aim we use pointwise estimates of the velocity and of  the pressure.

{ Mathematical analysis of the problem of the Navier-Stokes equations with artificial boundary condition was performed by many authors but without considering the rotation of body, see e.g.
\cite{BoFa,BrMu,DeKr,DeKr1}. The article can be seen as a first result in the case of motion of viscous fluids around rotating and translating body with artificial boundary condition.}

\medskip

The paper is organized as follows:  In the rest of this section we introduce notation and give some auxiliary results. The
next section 2 deals with pointwise estimates of the pressure of  the linear system (\ref{eq_1.1}). In Section 3  we consider the linear system     (\ref{eq_1.1}) with artificial boundary conditions. The error estimate of the velocity is derived comparing to the solution to the system given in  the exterior domain.
First let us introduce notation:\
\subsection*{Definitions and notation related to the rotational system}

Define
$  s(y):=1+|y|-y_1$ for $y \in \mathbb{R}^3 $,

\medskip

$\mathcal D_R:=B_{R}\setminus\overline {\mathcal D}$,

\medskip
 $B_R^c:=\mathbb R^3\setminus\overline{B_R}$,

\medskip\noindent
where $B_R:=\{x\in\mathbb R^3 ; |x|<R\},$\ for $R>0$ such that $B_R \supset\overline{\mathcal D.}   $

\noindent
 So, $\mathcal D_R$ is the truncation  of the exterior domain $ \overline{\mathcal D}^c:=\mathbb R ^3\setminus\overline{\mathcal D}\, $ by the ball $B_R$. The boundary $\partial \mathcal D_R$ consists of parts $\partial \mathcal D$ and $\partial B_R$, the later we  call the truncating  boundary.

\medskip

Fix $\tau \in (0,\infty)$, $\rho \in \mathbb R \setminus \{0\}$,
and put
$e_1:=(1,0,0),\;\Omega : = \rho \begin{pmatrix} 0&0&0\\ 0&0&-1\\ 0&1&0 \end{pmatrix}$,

so that
$\color{black}
\Omega \cdot  z = \rho e_1 \times z
\color{black}\
$
for $z\in \mathbb R^3$.

For $U \subset\mathbb R^3$ open, $u\in W^{2,1}_{{\rm loc}}(U)^3$, $z\in U$, put
\begin{align*}
(Lu) (z):= & -\Delta u(z) + \tau \partial_1 u(z) - (\rho e_1 \times z) \cdot \nabla u (z)
+\rho e_1 \times u(z),\qquad \qquad \  \\
\ \ \ (L^*u)(z):= & - \Delta u(z) - \tau \partial_1 u(z) + (\rho e_1 \times z) \cdot \nabla u(z)
-\rho e_1 \times u(z).
\end{align*}

Put
\begin{eqnarray*}
&&K(z,t):= (4\pi t)^{-3/2}e^{-|z|^2/(4t)} \quad (z\in\mathbb R^3, t\in (0,\infty)),\\
&&\Lambda (z,t):=\left(K(z,t)\delta_{jk}+ \partial z_j \partial z_k
\left(\int_{\mathbb R^3}(4\pi |z-y|)^{-1} K(y,t) dy\right)\right)_{1\le j,k\le 3} \qquad \\
&&(z\in \mathbb R^3,\ t>0),\\
&&
\Gamma(x,y,t):=
\color{black}
\Lambda (x-\tau  te_1 - e^{-t\Omega}y,t)
\color{black}
\cdot e^{-t\Omega},\\
&&
\widetilde{\Gamma}(x,y,t):= \Lambda (x+\tau  te_1 - e^{t\Omega}y,t)\cdot e^{t\Omega}
\quad (x,y\in\mathbb R^3, t>0),\\
&&\mathcal Z(x,y):= \int^{\infty}_0 \Gamma (x,y,t) dt,
\ \widetilde{ \mathcal Z}(x,y):=\int^{\infty}_0 \widetilde{\Gamma} (x,y,t) dt,
\\
&&
%\ \
(x,y\in\mathbb R^3, \, x\not = y).
\end{eqnarray*}
For $q\in (1,2)$, $f\in L^q(\mathbb R^3)^3$, put
$$
\mathcal R(f)(x) : = \int_{\mathbb R^3} \mathcal Z(x,y) f(y) dy \quad (x\in \mathbb R^3);
$$
see \cite[Lemma 3.1]{2}.
We will use the space

\medskip
$D^{1,2}_0( \overline{ \mathcal D}^c)^3:=\{v \in L^6( \overline{ \mathcal D}^c)^3
\cap H^1_{loc}( \overline{ \mathcal D}^c)^3\, :\,
\nabla v \in L^2( \overline{ \mathcal D}^c)^9,\; v| \partial \mathcal D =0\}$ equipped with the norm
$\|\nabla u\|_2$, where $v|\partial\mathcal D$ means the trace of $v$ on $\partial \mathcal D$.

For $p \in (1, \infty ),$ define $M _p$ as the space of all pairs
of functions $(u,\pi )$ such that
$
u \in W^{2,p}_{loc}( \overline{ \mbox{$\mathcal D$} }^c)^3,\;
\pi \in W^{1,p}_{loc}( \overline{ \mbox{$\mathcal D$} }^c),
$
\begin{eqnarray*} &&
u| \mbox{$\mathcal D$} _R \in W^{1,p}( \mbox{$\mathcal D$} _R)^3,
\quad
\pi | \mbox{$\mathcal D$} _R \in L^p( \mbox{$\mathcal D$} _R),
\quad
u| \partial \mbox{$\mathcal D$} \in W^{2-1/p,\,p}( \partial \mbox{$\mathcal D$} )^3,
\\ && \hspace{3em} \nonumber
\mbox{div}\, u| \mbox{$\mathcal D$} _R \in W^{1,p}( \mbox{$\mathcal D$} _R),
\quad
L(u)+ \nabla \pi | \mbox{$\mathcal D$} _R \in L^p( \mbox{$\mathcal D$} _R)^3
\end{eqnarray*}
for some $ R \in (0, \infty ) $ with $ \overline{ \mbox{$\mathcal D$} } \subset B_R$.

We write $C$ for generic constants.
It should be clear from context which are the parameters these constants depend on.
In order to lift possible ambiguities,
we sometimes use the notation $C (\gamma _1,\, ...,\, \gamma _n)$
in order to indicate that the constant in question depends in particular on
$\gamma _1,\, ...,\, \gamma _n \in (0, \infty ) $, for some $n \in \mathbb{N} $.
But the relevant constant may depend on other parameters as well.

\subsection*{Auxiliary results to  asymptotic behavior of the pressure}

\begin{lem}\label{Theorem 1} {\rm (Weyl's lemma).}
%\it
Let $n\in\mathbb N$,  $U\subset \mathbb R^n$ open, $u \in L^1_{{\rm loc}}(U)$ with $\int_Uu\cdot \Delta l\, dx =0$ for $l\in C^{\infty}_0(U)$. Then $u\in C^{\infty}(U)$ and $\Delta u =0$.

\end{lem}
{\it Proof:} An elementary proof is given in
\cite[Appendix]{4} \hfill$\Box$

For $q\in (1,3/2)$, $h\in L^q(\mathbb R^3)$, put
\[
\mathcal N(h)(x) := \int_{\mathbb R^3} - (4\pi|x-y|)^{-1}h(y) dy \quad (x\in \mathbb R^3).
\]
For $q\in (1,3)$, $h\in L^q(\mathbb R^3)$, put
\[
\mathcal S(h)(x):=\left(\int_{\mathbb R^3}(4\pi|x-y|^3)^{-1} (x-y)_j\cdot h(y)dy\right)_{1\le j\le 3}
\quad (x\in\mathbb R^3).
\]
For $q\in (1,3)$, $h\in L^q(\mathbb R^3)^3$, put
\[
\mathcal P(h)(x):=\int_{\mathbb R^3} (4\pi|x-y|^3)^{-1} ((x-y)\cdot h(y))dy \quad (x\in\mathbb R^3).
\]

Note that $S(h)$ is a vector-valued function with $h$ being scalar, whereas $P(h)$ is a scalar
function with $h$ being vector-valued.

\begin{lem}\label{Theorem 2}
{\it Let $q\in (1,3/2)$, $h\in L^q(\mathbb R^3)$.
Then $\mathcal N(h) \in W^{2,q}_{{\rm loc}} (\mathbb R^3) \cap L^{(1/q-2/3)^{-1}}(\mathbb R^3)$,
$\Delta\mathcal N(h)=h$.
If $h\in W^{1,q}(\mathbb R^3)$,
then $\partial_l \mathcal N(h) = \mathcal N(\partial_lh)$ $(1\le l \le 3)$.

Let $q\in (1,3)$, $h\in L^q (\mathbb R^3)$. Then
$
\mathcal S(h) \in W^{1,q}_{{\rm loc}}(\mathbb R^3)^3
,
\ {\rm div}\,\mathcal S(h)=h.
$
If $q\in (1,3/2)$,
then $\nabla \mathcal N(h) = \mathcal S(h)$.
If $h\in W^{1,q}(\mathbb R^3)$,
then $\mathcal S(h)\in W^{2,q}_{{\rm loc}}(\mathbb R^3)^3
$.

Let $q\in (1,3)$, $h\in L^q(\mathbb R^3)^3$. Then
\begin{align*}
&\mathcal P(h)\in W^{1,q}_{{\rm loc}}(\mathbb R^3) \cap L^{(1/q-1/3)^{-1}}(\mathbb R^3),
\\&
\left(\int_{\mathbb R^3}\left(\int_{\mathbb R^3}|x-y|^{-2}|h(y)|dy\right)^{(1/q-1/3)^{-1}}dx\right)^{1/q-1/3}
\le C\|h\|_q.
\end{align*}
}
\end{lem}
{\it Proof:} The assertion of the Lemma \ref{Theorem 2} follows from well known Hardy-Littlewood-Sobolev inequality, Calderon-Zyg\-mund inequality,
and density arguments.\hfill $\Box$

\begin{lem}\label{Fa}
\cite[Lemma 2.2]{KrNoPo}
Let $B\in\mathbb R,\, S\in(0,\infty).$
Then
\begin {equation}\displaystyle \int
_{\partial B_R}s(x)^{-B}
\, do_x\le C(S,B)\cdot R^{2-\hbox{min}\,\{1,B\}}\cdot\sigma
(R)
\end{equation}
for $R\in[S,\infty),$ with $\sigma(R):=1$
if $B\ne1,$ and $\sigma(R)={\rm
ln(1+R)}$ if $B=1$.

\end{lem}

\section{Decay estimates }

In first part of this section we  recall some known results from \cite{2} and \cite{3}
about the decay of the velocity part of the solution of the system (\ref{eq_1.1}), and in order to get the full decay characterization of the solution we derive the decay of the pressure part of solution of (\ref{eq_1.1}). In the second part of this section we extend the result for the pressure to the non-linear case of (\ref{1.0}).

\medskip

\subsection*{Decay estimates in the linear case}

Our starting point is a decay result from \cite{3} for the velocity part $u$ of a solution to
(\ref{eq_1.1}).

\begin{thm}\label{thm_1.1}
{\rm (\cite[Theorem 3.12]{3})}
Suppose that \mbox{$\mathcal D$} is $C^2$-bounded.
Let $p\in (1,\infty)$, $(u,\pi)\in M_p$. Put $F = L(u)+\nabla \pi$.
Suppose there are numbers $S_1,S,\gamma \in (0,\infty)$, $A\in [2,\infty)$, $B\in \mathbb R$
  such that $S_1 < S$,
  \[\overline{\mathcal D}\cup {\rm supp}({\rm div}\,u) \subset B_{S_1}, \quad u|{B_{S}^c} \in L^6 (B^c_S)^3,
  \quad
  \nabla u|{B^c_S} \in L^2(B_S^c)^9,
  \]
\[
A+\min \{1,B\} \ge 3, \ |F(z)|\le \gamma|z|^{-A}s(z)^{-B} \ \text{for} \ z\in B^c_{S_1}.
\]
Then
\begin{equation}\label{1.6}
\displaystyle|u(y)|\le
C\, (|y|s(y))^{-1}\,l_{A,B}(y),
\end{equation}
\begin{equation}\label{1.7}
| \nabla  u(y)| \le
C\, (|y|s(y))^{-3/2} \, s(y)^{ \max{(0,7/2-A-B)}}
\color{black}
\,l_{A,B}(y)
\color{black}
\end{equation}
for $y\in B^c_S$,
where
function $l_{A,B}$
is given by
$$
\left\{
\begin{array}{crl}
1 \quad & \text{if} \quad A + \min \{1,B\} > 3\\
\max (1, \text{ln}(y)) & \text{if} \quad A + \min \{1,B\} = 3.
\end{array} \right.
$$
\end{thm}
The requirements $u|{B_{S}^c} \in L^6 (B^c_S)^3,\; \nabla u|{B^c_S} \in L^2(B_S^c)^9$
should be interpreted as decay conditions on $u$.

It may be deduced from Theorem \ref{thm_1.1} that inequalities
(\ref{1.6}) and (\ref{1.7}) hold under assumptions weaker than those stated in that
theorem. We specify this more general situation in the ensuing corollary, which in addition
indicates some properties of $F$ that will be useful in the following.

\begin{cor} \label{corollary2.1}
  Let $p \in (1, \infty ),\; \gamma ,\,S_1,\, S \in (0, \infty ) $ with
  $\overline{ \mbox{$\mathcal D$} }\subset B_{S_1},\; S_1<S,\; A \in [2, \infty ),\;B \in \mathbb{R} $
    with $A+\min\{1,B\}\ge 3$. Let $F:\overline{ \mbox{$\mathcal D$} }^c \mapsto \mathbb{R}^3 $
    be measurable with $F| \mbox{$\mathcal D$} _{S_1} \in L^p( \mbox{$\mathcal D$} _{S_1})^3$
    and
    $|F(z)|\le \gamma|z|^{-A}s(z)^{-B} \ \text{for} \ z\in B^c_{S_1}.$

    Let $u \in W^{1,p}_{loc}( \overline{ \mbox{$\mathcal D$} }^c)^3$ with
    $u|{B_{S}^c} \in L^6 (B^c_S)^3$,
$\nabla u|{B^c_S} \in L^2(B_S^c)^9$,
    ${\rm supp}({\rm div}\,u)\subset B_{S_1}$,
    \begin{eqnarray} \label{2.*}&&\hspace{-3em}
      \int_{ \overline{ \mbox{$\mathcal D$} }^c}\bigl[\, \nabla u \cdot \nabla \varphi
        + \bigl(\, \tau \, \partial _1u  - ( \varrho \, e_1 \times z) \cdot \nabla u
        + ( \varrho \, e_1 \times u) - F \,\bigr) \cdot \varphi \,\bigr]\, dz
      \\&&\nonumber \hspace{-3em}
      =0 \quad \mbox{for}\;\;
    \varphi  \in C ^{ \infty } _0( \overline{ \mbox{$\mathcal D$} }^c)^3\;\;\mbox{with}\;\;
    {\rm div}\, \varphi =0.
    \end{eqnarray}
Then inequalities (\ref{1.6}) and (\ref{1.7}) hold for $y \in B_S^c$.

Moreover $F \in L^q(\overline{ \mbox{$\mathcal D$} }^c)^3$ for $q \in (1,p]$.
  If $p\ge 6/5,$ the function $F$ may be considered as a bounded linear functional on
  $\mbox{$\mathcal D$} ^ {1,2}_0( \overline{ \mbox{$\mathcal D$} }^c)^3$, in the usual sense.

\color{black}
Let $\pi \in L^p_{loc}( \overline{ \mbox{$\mathcal D$}}^c)$ with
\begin{eqnarray} \label{2.**}&&\hspace{-2em}
      \int_{ \overline{ \mbox{$\mathcal D$} }^c}\bigl[\, \nabla u \cdot \nabla \varphi
        + \bigl(\, \tau \, \partial _1u  - ( \varrho \, e_1 \times z) \cdot \nabla u
        + ( \varrho \, e_1 \times u) - F \,\bigr) \cdot \varphi
        \\&&\nonumber
        -\pi\, {\rm div}\, \varphi \,\bigr]\, dz
      =0 \quad \mbox{for}\;\;
    \varphi  \in C ^{ \infty } _0( \overline{ \mbox{$\mathcal D$} }^c)^3.
    \end{eqnarray}
Fix some number $S_0 \in (0,S_1)$ with $\overline{ \mbox{$\mathcal D$}}\cup {\rm supp}({\rm div}\, u)
\subset B_{S_0}.$
Then the relations
$u| \overline{ B_{S_0}}^c \in W^{2,p}_{loc}(\overline{B_{S_0}}^c)^3,\;
\pi \in W^{1,p}_{loc}(\overline{B_{S_0}}^c)$
and
$L(u| \overline{ B_{S_0}}^c)+\nabla \pi = F| \overline{ B_{S_0}}^c$
hold.
\color{black}
  \end{cor}
  {\it Proof:}
  For $z \in B_{S_1}^c$, we have
  \begin{eqnarray*}&&\hspace{-1em}
    |F(z)|
    \le
    \gamma \, C(S_1,A)\,|z| ^{-2} \, s(z)^{-A+2-B}
    \le
    \gamma \,C(S_1,A)\,|z| ^{-2} \, s(z) ^{-A+2-\min\{1,B\}}
    \\[1ex]&&\hspace{-1em}
    \le
    \gamma \,C(S_1,A)\, |z| ^{-2} \, s(z) ^{-1}.
    \end{eqnarray*}
  Thus for $q \in (1, \infty ),$ with
  Lemma \ref{Fa},
  \begin{eqnarray*}
    \int_{ B_{S_1}^c}|F(z)|^q\, dz
    \le
    C\, \int_{ S_1} ^{ \infty } r^{-2q}\, \int_{ \partial B_r}s(z) ^{-q}\, do_z\, dr
    \le
    C\, \int_{ S_1} ^{ \infty } r^{-2q+1}\,dr < \infty .
  \end{eqnarray*}
  It follows that $F \in L^q( \overline{ \mbox{$\mathcal D$}}^c)^3$ for $q \in (1,p]$.
According to \cite[Theorem II.6.1]{Galdi}, the inequality $\|v\|_6\le C\, \| \nabla v\|_2$
holds for $v \in
\color{black}
D ^{1,2}_0(\overline{ \mbox{$\mathcal D$}}^c)^3
\color{black}
$.
Thus, if $p\ge 6/5$, hence $F \in L^{6/5}(\overline{ \mbox{$\mathcal D$}}^c)^3$,
this function $F$  may be considered as a linear bounded functional on
$\mbox{$\mathcal D$} ^{1,2}_0(\overline{ \mbox{$\mathcal D$}}^c)^3$.
The $L^p$-integrability of $F$ and the assumptions on $u$ imply that the function
\begin{eqnarray} \label{BB}
  G(z):=F(z)-
\bigl(\, \tau \, e_1  - ( \varrho \, e_1 \times z) \,\bigr)  \cdot \nabla u(z)
- \bigl(\, \varrho \, e_1 \times u(z) \,\bigr) ,\;\;
z \in \overline{ \mbox{$\mathcal D$}}^c,
\end{eqnarray}
belongs to $L^p_{loc}(\overline{ \mbox{$\mathcal D$}}^c)^3$.
\color{black}
The choice of $S_0$ (see at the end of Corollary \ref{corollary2.1}) means in particular that
${\rm div}\,(u| \overline{ B_{S_0}}^c)=0.$
This equation,
(\ref{2.**}),
the relation $G \in L^p_{loc}(\overline{ \mbox{$\mathcal D$}}^c)^3$
and interior regularity of solutions to the Stokes system (see \cite[Theorem IV.4.1]{Galdi}
for example) imply
the claims in the last sentence of Corollary \ref{corollary2.1}.
\color{black}

%that
%$u| \overline{ B_{S_0}}^c \in W^{2,p}_{loc}(\overline{B_{S_0}}^c)^3$ and there is
%$\pi \in W^{1,p}_{loc}(\overline{B_{S_0}}^c)$
%with
%$L(u| \overline{ B_{S_0}}^c)+\nabla \pi = F| \overline{ B_{S_0}}^c.$
Put $ S_0 ^{\prime} :=(S_0+S_1)/2,\;
A_{S_0 ^{\prime} ,R}:=B_R \backslash B_{S_0 ^{\prime} } $ for $R \in (S_0 ^{\prime} , \infty ).$
Then
$u| A_{S_0 ^{\prime} ,R} \in W^{2,p}(A_{S_0 ^{\prime} ,R})^3$ and
$\pi | A_{S_0 ^{\prime} ,R} \in W^{1,p}(A_{S_0 ^{\prime} ,R})^3$
for $R \in (S_0 ^{\prime} , \infty ),$
so
$(u| A_{S_0 ^{\prime} ,R},\, \pi| A_{S_0 ^{\prime} ,R}) \in M_p$,
with $B_{S_0 ^{\prime} } $ in the role of \mbox{$\mathcal D$}.
Note that $S_0<S_o ^{\prime} <S_1< S.$ Thus the assumptions of Theorem \ref{thm_1.1} are
satisfied with \mbox{$\mathcal D$}  replaced by $B_{S_0 ^{\prime} }$. As a consequence inequalities
(\ref{1.6}) and (\ref{1.7}) hold.
\hfill $\Box $
\noindent
\begin{rem}\label{rem_existence}  Solutions as considered in
Corollary \ref{corollary2.1}
  exist if, for example,
Dirichlet boundary conditions are prescribed on $\partial \mathcal D$. In fact,
as stated in \cite[Theorem VIII.1.2]{Galdi},
if $F$ is a bounded linear functional on the space
$D^{1,2}_0( \overline{ \mathcal D}^c)^3$,
and if $b \in H ^{1/2} ( \partial \overline{ \mathcal D}^c)^3,$
then
there is a function
$u \in L^6( \overline{ \mathcal D}^c)^3\cap  W^{1,1}_{loc}( \overline{ \mathcal D}^c)^3
$
such that
$ \nabla u \in L^2( \overline{ \mathcal D}^c)^9$
and
$u$ satisfies the equations (\ref{2.*}) and ${\rm div}\,u =0$ (weak form of
(\ref{eq_1.1})), as well as the boundary conditions $u| \partial \mbox{$\mathcal D$} =b$.
\color{black}
Existence of a pressure $\pi \in L^p_{loc}( \overline{ \mathcal D}^c)$ with (\ref{2.**})
holds according to \cite[Lemma VIII.1.1]{Galdi}.
\color{black}
\end {rem}
\medskip
The main result of this section, dealing with the asymptotics of the pressure, is stated in
\begin{thm}\label{thm_Main_1}
{\it
Let $p,\, \gamma ,\, S_1,\, S,\,  A,\, B,\, F,\, u$
be given as in Corollary \ref{corollary2.1}, but with the stronger assumptions
$A=5/2,\; B \in (1/2, \, \infty )$ on $A$ and $B$.
\color{black}
Let $\pi \in L^p_{loc}( \overline{ \mbox{$\mathcal D$}}^c)$
such that (\ref{2.**}) holds
\color{black}
Then there is $c_0 \in \mathbb{R} $ such that
\begin{eqnarray} \label{CC}
|\pi(x)+c_0|\le C \, |x| ^{-2} \quad \mbox{for}\;\; x \in
B_S^c.
\end{eqnarray}
}
\end{thm}
{\it Proof:}
\color{black}
By Corollary \ref{corollary2.1} we have $F \in L^q(\overline{ \mathcal D}^c)^3$ for $q \in (1,p]$.
Fix some number $S_0 \in (0,S_1) $ with
$\overline{ \mbox{$\mathcal D$}}\cup {\rm supp}({\rm div}\, u)\subset B_{S_0}$.
Then again by Corollary \ref{corollary2.1}, the relations
$u| \overline{ B_{S_0}}^c \in W^{2,p}_{loc}( \overline{ B_{S_0}}^c)^3,
\;\pi | \overline{ B_{S_0}}^c \in W^{1,p}_{loc}( \overline{ B_{S_0}}^c)$
and $L(u| \overline{ B_{S_0}}^c )+\nabla (\pi | \overline{ B_{S_0}}^c )=F | \overline{ B_{S_0}}^c $ hold.
\color{black}
Note that $S_0<S_1<S$.
Take $\phi \in C^{\infty}(\mathbb R^3)$ with
$$
%\phi|B_{S+\frac14(R_0-S)}=0, \ \phi|B^c_{S+\frac34(R_0-S)} = 1,
\phi|B_{S_1+\frac14(S-S_1)}=0, \ \phi|B^c_{S_1+\frac34(S-S_1)} = 1,
$$
and put
$\widetilde u : = \phi  \cdot u$, $\widetilde{\pi} := \phi \cdot \pi$,
with $\widetilde u, \widetilde{\pi}$ to be considered as functions in $\mathbb R^3$.
By the choice of $\phi $ and the properties of $u$ and $\pi$,
we get
$\widetilde u \in W^{2,q}_{{\rm loc}}(\mathbb R^3)^3$,
$\widetilde{\pi}\in W^{1,q}_{{\rm loc}}(\mathbb R^3)$ for $q\in [1,p]$,
$\widetilde u |B_S^c=u|B_S^c\in L^6(\mathbb R^3)^3$,
$\nabla \widetilde u|B_S^c=\nabla u|B_S^c \in L^2(\mathbb R^3)^9$. Put
\begin{align*}
g_l(z):= &
\color{black}
- \sum^3_{k=1}  \partial_k \phi(z) \partial_k u_l(z) - \Delta \phi(z) u_l(z)
+ \tau \partial_1 \phi(z) u_l(z)
\color{black}
\\
&- \sum^3_{k=1} (\tau e_1 \times z)_k \cdot \partial_k \phi(z) \cdot u_l (z) + \partial_l \phi(z) \pi(z)
\end{align*}
for $z\in\mathbb R^3$, $1\le l \le 3$, and set $\gamma :={\rm div}\, \widetilde{ u}$. Then
\begin{align}
&{\rm supp}(g) \subset \overline{ B_{S_1+3(S-S_1)/4}} \setminus B_{S_1+(S-S_1)/4},
\ g\in L^q(\mathbb R^3)^3 \ \text{for} \ q\in [1,p], \nonumber
\\&
L\widetilde u + \nabla \widetilde{\pi} = g + \phi \cdot F,
\
\gamma
= \nabla \phi \cdot u, \label{eq1}
\end{align}
in particular $\mbox { supp}
( \gamma )
\subset \overline{ B_{S_1+3(S-S_1)/4}} \setminus B_{S_1+(S-S_1)/4}$,
$
\color{black}
\gamma
\in W^{2,q}(\mathbb R^3),
\color{black}
$
$g+ \phi \cdot F \in L^q(\mathbb R^3)^3$ for $q\in (1,p]$.
Let $x\in \mathbb R^3$, $\varepsilon >0$ with $\overline{ B_{\varepsilon}(x)}\subset B_{S_1}$, where $B_\varepsilon(x)=\{y\in\mathbb R^3;|y-x|<\varepsilon\}$.
Since $\widetilde u|B_{S_1} = 0$, $\widetilde{\pi}|B_{S_1}=0$,
it follows from \cite[Theorem 3.11]{3} with $\mathcal D$ replaced by $B_{\varepsilon}(x)$  that
$$
\widetilde u(y) = \mathcal R(g+\phi F)(y) + \mathcal S(
\gamma
)(y)
\ \text{for} \ y\in \overline{B_{\varepsilon}(x)}^c.
$$
Since this is true for any $x\in \mathbb R^3$, $\varepsilon >0$ with
$\overline{B_{\varepsilon}(x)} \subset B_{S_1}$, it follows that
\begin{equation}\label{eq2}
\widetilde u = \mathcal R(g+\phi F)+\mathcal S(
\gamma
) \ \text{in} \ \mathbb R^3.
\end{equation}
But $\mathcal S(
\gamma
) \in W^{2,q}_{{\rm loc}}(\mathbb R^3)^3$
for $q\in [1,\min \{3,p\})$ by Lemma \ref{Theorem 2}, so from \eqref{eq2}
$$
\mathcal R(g+ \phi F)\in W_{{\rm loc}}^{2,q}(\mathbb R^3)^3 \ \text{for} \ q \in [1,\min\{3,p\}).
$$
This relation and
\cite[(3.11) and the
\color{black}
inequality
\color{black}
following (3.15)]{2}
imply
\begin{equation}\label{eq3}
\sup_{z\in\mathbb R^3} \int_{\mathbb R^3} |\mathcal Z(z,y)\cdot (g+ \phi F)(y)|dy<\infty.
\end{equation}
Let $\psi\in C^{\infty}_0 (\mathbb R^3)^3$. Due to \eqref{eq3}, we may apply Fubini's theorem, to obtain
\begin{align}\label{eq4}
A&:=\int_{\mathbb R^3}\psi(x) (L\mathcal R(g+\phi F))(x)dx
= \int_{\mathbb R^3}(L^*\psi)(x) \mathcal R(g+\phi F)(x)dx
\\
&=\int_{\mathbb R^3}\int_{\mathbb R^3}[(L^*\psi)(x)]^T\cdot \mathcal Z(x,y) \cdot (g+ \phi F)(y) dydx\nonumber\\
&= \int_{\mathbb R^3}\int_{\mathbb R^3} [(L^*\psi)(z)]^T \cdot \mathcal Z(x,y) \cdot (g+\phi F)(y) dx dy.\nonumber
\end{align}
But for $a,b,x,y \in\mathbb R^3$ with $x\not=y$,
$$
a^T\cdot \mathcal Z(x,y) \cdot b = \int^{\infty}_0 a^T \Gamma (x,y,t) b dt,
$$
hence with \cite[Lemma 2.10]{2},
\begin{align*}
a^T\cdot \mathcal Z(x,y) \cdot b &= \int^{\infty}_0 a^T \cdot e^{-t\Omega} \Lambda (e^{t\Omega}x - \tau te_1 - y,t) \cdot b dt\\
&= \int^{\infty}_0 b^T [e^{-t\Omega} \Lambda (e^{t\Omega} x - \tau t e_1 - y,t) ]^T a dt\\
&= \int^{\infty}_0 b^T \cdot \Lambda (e^{t\Omega} x - \tau te_1 - y,t) e^{t\Omega} \cdot a dt\\
&= \int^{\infty}_0 b^T \Lambda (y + \tau te_1 - e^{t\Omega} x,t) e^{t\Omega} a dt\\
&= \int^{\infty}_0 b^T \widetilde {\Gamma} (y,x,t) \cdot a dt = b^T \cdot \widetilde Z(y,x)\cdot a.
\end{align*}
Therefore from \eqref{eq4}
\begin{equation}\label{eq5}
A = \int_{\mathbb R^3} (g+ \phi F)(y)^T\cdot \int_{\mathbb R^3} \widetilde{\mathcal Z}(y,x) \cdot (L^*\psi)(x)dxdy.
\end{equation}
Since $\psi\in C^{\infty}_0(\mathbb R^3)^3$, we may choose $x_0 \in \mathbb R^3$, $\varepsilon >0$ such that
$$
\overline{B_{\varepsilon}(x_0)}\subset \mathbb R^3 \setminus {\rm supp}\,(\psi).
$$
Thus we get from \cite[Theorem 4.3]{1}
\color{black}
with $\mathcal D,U,\omega,u,L$ replaced by $B_{\varepsilon}(x_0)$,
$\tau e_1$,
$-\rho e_1,\, \psi,\, L^*$, respectively,
\color{black}
and with $\pi=0$, that
$$
\int_{\mathbb R^3} \widetilde{\mathcal Z}(y,x)\cdot (L^*\psi)(x) dx = \psi(y) - \mathcal S ({\rm div}\,\psi)(y)
$$
for $y\in\mathbb R^3 \setminus\overline{B_{\varepsilon}(x_0)}$. Since this is true for any
$x_0 \in\mathbb R^3$,
$\varepsilon >0$ with $\overline{B_{\varepsilon}(x_0)} \subset \mathbb R^3 \setminus {\rm supp}(\psi)$,
the preceding equation holds for any $y\in \mathbb R^3$. It follows from~\eqref{eq5}
\begin{equation}\label{eq6}
\int_{\mathbb R^3}\psi(x)(L\mathcal R(g+ \phi F))(x) dx
=\int_{\mathbb R^3} (g+ \phi F)(y) \cdot (\psi(y) - \mathcal S({\rm div}\,\psi)(y)) dy.
\end{equation}
Again recalling that $g+\phi F\in L^q(\mathbb R^3)^3$ for $q\in (1,p]$, we get with Lemma \ref{Theorem 2} that
\begin{equation}\label{eq7}
\int_{\mathbb R^3}\psi(x) \nabla \mathcal P(g+\phi F)(x)dx
= \int_{\mathbb R^3} - {\rm div}\,\psi(x) \cdot \mathcal P(g+\phi F)(x)dx.
\end{equation}
Put $q_0 := \min \{6/5,p\}$, and note that $q_0\in (1,3/2)$, $q_0 \le p$.

Thus, by H\"older's inequality
and Lemma \ref{Theorem 2},
\begin{align*}
&
\int_{\mathbb R^3}|{\rm div}\,\psi(x)\, \mathcal P(g+
\color{black}
\phi
\color{black}
\, F)(x)|\, dx
\\
&\le \int_{\mathbb R^3}\int_{\mathbb R^3} |{\rm div}\,\psi(x)
(4\pi|x-y|^3)^{-1}  (x-y)\cdot (g+\phi F)(y)| dydx
\\[1ex]
&\le \|{\rm div}\,\psi\|_{(4/3-1/q_0)^{-1}}
\\&\hspace{2em}
\left(\int_{\mathbb R^3}
\left(\int_{\mathbb R^3}(4\pi|x-y|^2)^{-1} |(g+\phi F)(y)| dy\right)^{(1/q_0-1/3)^{-1}}dx\right)^{1/q_0-1/3}\\
&\le C \cdot \|{\rm div}\,\psi\|_{(4/3-1/q_0)^{-1}}\cdot \|g+\phi F\|_{q_0} < \infty.
\end{align*}
As a consequence, we may apply Fubini's theorem to deduce from (\ref{eq7}) that
\begin{align}
&\int_{\mathbb R^3} \psi(x)\nabla \mathcal P(g+\phi F)(x)dx\label{eq8}
\\
&=-\int_{\mathbb R^3}\int_{\mathbb R^3} ({\rm div}\,\psi)(x)
(4\pi|x-y|^{3})^{-1}(x-y)\cdot (g+\phi F)(y)dxdy\nonumber\\
&=\int_{\mathbb R^3} (g+\phi F)(y) \cdot \mathcal S({\rm div}\,\psi)(y)dy.\nonumber
\end{align}
From \eqref{eq6} and \eqref{eq8},
\begin{align*}
&\int_{\mathbb R^3} \psi(x) ((L\mathcal R(g+\phi F))(x) + \nabla \mathcal P(g+\phi F)(x))dx\\
&= \int_{\mathbb R^3} \psi(x) (g+\phi F)(x) dx.
\end{align*}
Since this is true for any $\psi\in C^{\infty}_0 (\mathbb R^3)^3$, we have found that
\begin{equation}\label{eq9}
L\mathcal R(g+\phi F)+\nabla \mathcal P(g+\phi F)=g+\phi F.
\end{equation}
On the other hand, by \eqref{eq2} and \eqref{eq1}
$$
L\mathcal R(g+\phi F)+L\mathcal S(
\gamma
) + \nabla \widetilde{\pi} = L\widetilde u + \nabla \widetilde {\pi} = g+ \phi F.
$$
By subtracting this equation from \eqref{eq9}, we get
\begin{equation}\label{eq10}
\nabla\mathcal P(g+\phi F)-L\mathcal S(
\gamma
) - \nabla \widetilde{\pi}=0.
\end{equation}
Next we consider the term
${\rm div} \bigl(\, L\mathcal S(\gamma ) \,\bigr) $.
Recall that
\color{black}
$q_0< 3/2,\; q_0\le p$ and
$ \gamma \in W^{2,q}(\mathbb R^3)$ for $q\in [1,p]$
(see (\ref{eq1})),
\color{black}
so by Lemma \ref{Theorem 2}
\begin{equation}\label{eq11}
\begin{cases}
\mathcal S(\gamma) \in W^{2,q_0}_{{\rm loc}}(\mathbb R^3)^3, \ {\rm div}\,\mathcal S(\gamma) = \gamma, \ \mathcal N(\gamma) \in W^{2,q_0}_{{\rm loc}} (\mathbb R^3),\\
\nabla \mathcal N(\gamma)=\mathcal S(\gamma).
\end{cases}
\end{equation}
\color{black}
Since
$e_1 \times \mathcal S( \gamma )=(0,\,-\mathcal S_3 ( \gamma ),\, \mathcal S_2( \gamma ))$
and because of the equation $\nabla \mathcal N( \gamma )=\mathcal S( \gamma )$
in (\ref{eq11}), we may conclude that
\begin{eqnarray} \label{2.17b}
{\rm div}(e_1 \times \mathcal S( \gamma ))
=
- \partial_2\mathcal S_3(\gamma) + \partial_3\mathcal S_2(\gamma)
=
- \partial_2\partial_3 \mathcal N(\gamma) + \partial_3\partial_2 \mathcal N(\gamma) = 0.
\end{eqnarray}
Moreover, for $z\in\mathbb R^3$, $1\le j\le 3$,
\begin{eqnarray} \label{2.17c}
(e_1\times z) \cdot \nabla {\mathcal S} _j(\gamma)(z)
=
-z_3\partial _2\mathcal S_j(\gamma)(z)+ z_2 \partial _3 \mathcal S _j(\gamma )(z).
\end{eqnarray}
Put
$\varphi (z):=-z_3 \partial_2 \gamma(z) + z_2\partial_3\gamma(z) $
for $z \in \mathbb{R}^3 .$
Then with (\ref{2.17c}), the equation ${\rm div}\mathcal S( \gamma )=\gamma $
in \eqref{eq11}, and the second and third equation in (\ref{2.17b}),
\begin{eqnarray*} &&\hspace{-2em}
{\rm div}_z \bigl(\, (e_1\times z)\cdot\nabla \mathcal S(\gamma)(z) \,\bigr)
=
\sum_{j = 1}^3 \partial z_j \bigl(\, (e_1\times z)\cdot\nabla \mathcal S_j(\gamma)(z) \,\bigr)
\\&&\hspace{-2em}
=
-z_3 \partial _2{\rm div}_z\mathcal S( \gamma )(z) +z_2 \partial _3{\rm div}_z\mathcal S( \gamma )(z)
+ \sum_{j = 1}^3 \bigl(\, \partial z_j(-z_3) \partial _2\mathcal S_j( \gamma )(z)
+\partial z_j(z_2) \partial _3\mathcal S_j( \gamma )(z) \,\bigr)
\\&&\hspace{-2em}
=
\varphi (z)- \partial _2\mathcal S_3( \gamma )(z)+\partial _3\mathcal S_2( \gamma )(z)
=
\varphi (z).
\end{eqnarray*}
\color{black}
Let $\psi\in C^{\infty}_0(\mathbb R^3)$. Then it follows that
\begin{align*}
&\int_{\mathbb R^3}\nabla \psi \cdot (\rho e_1 \times \mathcal S(\gamma)) dx =0,\\
&\int_{\mathbb R^3}\nabla \psi(z) [(\rho e_1\times z)\cdot \nabla \mathcal S(\gamma) (z)]dz
= \int_{\mathbb R^3} \psi(z)(-\varphi(z)) dz.
\end{align*}
Obviously, again with \eqref{eq11},
$$
\int_{\mathbb R^3}\nabla \psi \cdot \Delta \mathcal S(\gamma) dx = \int_{\mathbb R^3} \nabla \Delta \psi \cdot \mathcal S(\gamma) dx = -\int_{\mathbb R^3} \Delta\psi \cdot \gamma dx
= \int_{\mathbb R^3} \psi \cdot (-\Delta \gamma) dx,
$$
and similarly,
$$
\int_{\mathbb R^3} \nabla \psi(\tau \partial_1\mathcal S(\gamma))
\color{black}
dx
\color{black}
= \int_{\mathbb R^3} \psi(-\tau\partial_1 \gamma) dx.
$$
Combining these equations, we get
$$
\int_{\mathbb R^3} \nabla \psi \cdot L\mathcal S(\gamma) dx = \int_{\mathbb R^3} \psi(\varphi+ \Delta \gamma - \tau \partial_1 \gamma) dx.
$$
Now from \eqref{eq10}
\begin{equation}\label{eq12}
\int_{\mathbb R^3}\nabla \psi [\nabla \mathcal P(g+\phi F)-\nabla (\phi \pi)]dx
= \int_{\mathbb R^3}\psi(\varphi+ \Delta \gamma - \tau \partial_1 \gamma)dx.
\end{equation}
Since $\gamma \in W^{2,q}(\mathbb R^3)$ for $q\in [1,p]$
and ${\rm supp}(\gamma) \subset B_{S}\backslash B_{S_1}$
\color{black}
due to (\ref{eq1})),
\color{black}
it follows that
$\varphi+ \Delta \gamma - \tau \partial_1 \gamma \in L^q(\mathbb R^3)$ for $q\in [1,p]$q,
so we may consider
$\mathcal N (\varphi + \Delta \gamma - \tau \partial_1 \gamma)$. Lemma \ref{Theorem 2} yields
\begin{align*}
&\mathcal N(\varphi+ \Delta \gamma  - \tau \partial_1 \gamma) \in W^{2,q_0}_{{\rm loc}} (\mathbb R^3),\\
&\Delta \mathcal N(\varphi+ \Delta \gamma
- \tau\partial_1 \gamma) = \varphi+ \Delta \gamma - \tau \partial_1 \gamma.
\end{align*}
Therefore from \eqref{eq12}
$$
\int_{\mathbb R^3} \nabla \psi [\nabla \mathcal P(g+\phi F)- \nabla \mathcal N(\varphi+ \Delta \phi - \tau\partial_1 \gamma) - \nabla (\phi \pi) ] dx = 0.
$$
Lemma \ref{Theorem 1} now yields
\begin{eqnarray}
Q := \mathcal P(g+ \phi F)
- \mathcal N (\varphi+ \Delta \gamma - \tau \partial_1\gamma) - \phi \pi \in C^{\infty}(\mathbb R^3),
\quad
\label{99}
\Delta Q = 0.
\end{eqnarray}
Now we again apply Lemma \ref{Theorem 2}. Since $g+ \phi \cdot F \in L^{q_0}(\mathbb R^3)^3$, we have
$$
\mathcal P(g+ \phi F)\in L^{(1/q_0-1/3)^{-1}}(\mathbb R^3).
$$
Moreover $\varphi+ \Delta \gamma - \tau \partial_1 \gamma \in L^{q_0}(\mathbb R^3)$, so
$$
\mathcal N(\varphi+ \Delta \gamma - \tau \partial_1 \gamma) \in L^{(1/q_0 - 2/3)^{-1}}(\mathbb R^3).
$$
Since $q_0\le p,$ and in view of our remarks at the beginning of this proof we know that
$
u| \overline{ B_{S_0}}^c\in W^{2,q_0}_{loc}(  \overline{ B_{S_0}}^c),\;
\pi| \overline{ B_{S_0}}^c\in W^{1,q_0}_{loc}(  \overline{ B_{S_0}}^c),\;
L(u| \overline{ B_{S_0}}^c )+\nabla (\pi | \overline{ B_{S_0}}^c )=F | \overline{ B_{S_0}}^c ,\;
\mbox{div} ( u|  \overline{ B_{S_0}}^c) =0
$
and
$F\in L^{q_0}(\mathbb R^3)^3$.
By the choice of $u$ in Corollary \ref{corollary2.1},
we have
$u|B_{S}^c\in L^6(B_{S}^c)^3$.
\color{black}
Therefore
\color{black}
\cite[Theorem~2.1]{3} yields there is $c_0\in\mathbb R$ such that
$$
\pi + c_0|B^c_{2 S} \in L^{3q_0/(3-q_0)}(B_{2\,S}^c)+L^{3}(B^c_{2\, S}).
$$
But by (\ref{99}),
\begin{eqnarray*}
Q-c_0
=
\mathcal P(g+ \phi F)
- \mathcal N (\varphi+ \Delta \gamma - \tau \partial_1\gamma)
- \phi\,( \pi+c_0) + ( \phi -1)\, c_0,
\end{eqnarray*}
where ${\rm supp}( \phi -1) \subset B_{S}$ and ${\rm supp}( \phi ) \subset B_{S_1}^c$.
We may conclude that
\begin{eqnarray}
\label{100}&&
Q - c_0
\in L^{(1/q_0 - 1/3)^{-1}}( \mathbb{R}^3 ) + L^{(1/q_0-2/3)^{-1}}( \mathbb{R}^3 )
+ L^{q_0}(\mathbb R^3)
\\&&\nonumber \hspace{2em}
+ L^{3q_0/(3-q_0)}(\mathbb R^3) + L^{3}( \mathbb{R}^3 ).
\end{eqnarray}
Let $\varepsilon \in (0,\infty)$,
and let $(Q- c_0)_{\varepsilon}$
be the usual Friedrich's mollifier of $Q- c_0$ associated with $\varepsilon$.

Due to (\ref{99}), (\ref{100}) and by standard properties Friedrich's  mollifier,
the function $(Q- c_0)_{\varepsilon}$ is bounded and
$\Delta (Q- c_0)_{\varepsilon}=0$.
Now Liouville's theorem yields $(Q- c_0)_{\varepsilon}=0$.
Since this is true for any $\varepsilon >0$
\color{black}
and because $Q \in C ^{ \infty } ( \mathbb{R}^3 ),$
\color{black}
we may conclude that $Q- c_0=0$, that is,
\begin{eqnarray*}
\phi \,(\pi + c_0)=\mathcal P(g+\phi F)-\mathcal N(\varphi +\Delta \gamma - \tau \partial_1 \gamma)
+( \phi -1)\, c_0,
\end{eqnarray*}
hence
\begin{equation}\label{eq13}
\pi + c_0|B_{S}^c
=\mathcal P(g+\phi F)-\mathcal N(\varphi +\Delta \gamma - \tau \partial_1 \gamma)|B_{S}^c,
\end{equation}
where we used that ${\rm supp}( \phi -1) \subset B_{S}$ and $\phi |B_{S}^c=1$.
Since ${\rm supp}(g) \subset \overline{ B_{S_1+3(S-S_1)/4}}$, we have
\begin{equation}\label{eq14}
|\mathcal P(g)(x)| \le c \cdot |x|^{-2} \ \text{for} \ x \in B^c_{S}.
\end{equation}
Due to the assumptions
$A=5/2,\; B \in (1/2, \, \infty )$
and because $\phi  F|B_{S_1+(S-S_1)/4} = 0$ and $S_1<S_1+(S-S_1)/4<S$,
we get by \cite[Theorem 3.2]{Farwig2} or \cite[Theorem 3.4]{KrNoPo}
%\footnotemark[1] \footnotetext[1]
%{
%The
%result stated in both these theorems is not optimal,
%with respect of logarithmic terms:
%No logarithmic term appears  e.g. in the case $ A=5/2,B=1$.\ Nevertheless we need not this case.
%\   }
that
\begin{equation}\label{eq16}
|\mathcal P(\phi \, F)(x)| \le c \, |x|^{-2}\ \text{for} \ x\in B^c_{S}.
\end{equation}
\color{black}
Note that according to \cite[Theorem 3.4 (iii)]{KrNoPo}, a logarithmic factor
should be added on the right-hand
side of (\ref{eq16})
in the case $A=5/2,\; B=1$. But this factor is superfluous. In fact, if the relation
$|F(z)|\le \gamma |z|^{-A}s(z)^{-B}\;(z \in B_{S_1}^c)$
is valid with $A=5/2,\; B=1,$
it holds in the case $A=5/2,\; B=3/4,$ too. But then \cite[Theorem 3.4 (i), (iii)]{KrNoPo}
yields that (\ref{eq16}) holds as it is, without additional factor.
\color{black}

Define $\zeta (x)=-x_3  \gamma(x)$, $\tilde{\zeta }(x):=x_2\gamma (x)$ for $x\in\mathbb R^3$.
Then ${\rm supp}( \zeta ) \cup \,{\rm supp}(\tilde{\zeta})\subset B_{S}\setminus \overline{B_{S_1}}$,
\begin{eqnarray*}&&
\zeta,\tilde{\zeta} \in W^{2,q}(\mathbb R^3) \ \text{for} \ q\in [1,p],
\\&&
\varphi= \partial_2 \zeta + \partial_3 \tilde{\zeta}.
\end{eqnarray*}
It follows with Lemma \ref{Theorem 2} that
$$
\mathcal N(\varphi )
= \partial_2 \mathcal N(\zeta)+\partial_3\mathcal N(\tilde{\zeta})
= \mathcal S_2 (\zeta) + \mathcal S_3(\tilde{\zeta}).
$$
Similarly, since ${\rm supp}(\gamma)\subset B_{S}\setminus \overline{B_{S_1}}$,
$\gamma \in W^{2,q}(\mathbb R^3)$ for $q\in [1,p]$,
%\color{black}
$$
\mathcal N(\Delta \gamma -\tau \partial _1 \gamma )
= \sum^3_{k=1} \mathcal S_k(\partial_k\gamma)
- \tau \mathcal S_1(\gamma).
$$
Together
$$
\mathcal N(\varphi+ \Delta \gamma - \tau \partial_1 \gamma)
= \mathcal S_2 (\zeta) + \mathcal S_3(\tilde{\zeta})
+ \sum^3_{k=1} \mathcal S_k (\partial_k\gamma) - \tau \mathcal S_1 (\gamma).
$$
Since supp$(\zeta) \cup {\rm supp}(\tilde{\zeta}) \cup {\rm supp}(\gamma)
\subset B_{S_1+3(S-S_1)/4}$,
we may conclude that
\begin{equation}\label{eq17}
|\mathcal N(\varphi+ \Delta \gamma -\tau\partial_1 \gamma)(x)| \le C\,
|x|^{-2}\ \text{for} \ |\xi| \in  B^c_{S}.
\end{equation}
Inequality (\ref{CC}) follows from
\eqref{eq13}--\eqref{eq17}.%,
\hfill $\Box$

\vspace{1em}
We remark that Theorem \ref{thm_Main_1} remains valid if the assumptions on $A$ and $B$
are replaced by
the conditions $A\ge 5/2,\; A+\min\{1,B\}>1$, which are weaker than those in Corollary
\ref{corollary2.1}. This observation is made precise by the
ensuing corollary. Its proof is obvious, but this modified version of Theorem \ref{thm_Main_1}
still is interesting because its requirements on $A$ and $B$ are closer to the ones in
Corollary \ref{corollary2.1} than those stated in Theorem \ref{thm_Main_1}.
\begin{cor} \label{corollary2.2}
Let $p,\, \gamma ,\, S_1,\, S,\,  A,\, B,\, F,\, u$
be given as in Corollary \ref{corollary2.1}, but with the stronger assumptions
$A\ge 5/2,\; A+\min\{1,B\}>3$ on $A$ and $B$.
Let $\pi \in L^p_{loc}( \overline{ \mbox{$\mathcal D$}}^c)$ such that (\ref{2.**}) holds.
Then there is $c_0 \in \mathbb{R} $ such that inequality (\ref{CC}) is valid.
\end{cor}
{\it Proof:}
Put $B ^{\prime} := A-5/2+\min\{1,B\}$. Since $A+\min\{1,B\}>3$, we have $B ^{\prime} \in
(1/2,\, \infty ).$ Moreover, since $A\ge 5/2$, we find for
\color{black}
$z \in B_{S_1}^c$
\color{black}
that
\begin{eqnarray*}
|F(z)|
\le
\gamma \, C(S_1,A)\, |z|^{-5/2}\, s(z)^{-A+5/2-B}
\le
\gamma \,C(S_1,A)\, |z|^{-5/2}\, s(z)^{-B ^{\prime} }.
\end{eqnarray*}
Thus the assumptions of Theorem \ref{thm_Main_1} are satisfied with $B$ replaced by $B ^{\prime} $
and with a modified parameter $\gamma $. This implies the conclusion of Theorem \ref{thm_Main_1}.
\hfill $\Box $

\bigskip

\subsection*{Decay estimates in the non-linear case}

Let us assume now   the non-linear case, i.e. the system (\ref{1.0}).
First, recall the result about the decay properties of the velocity in this non-linear case:
\begin{thm}\label{nonlin}
\cite[Theorem 1.1]{DKN5}
Let $
\gamma ,\, S_1 \in (0, \infty ) ,\; p_0 \in (1, \infty ),\; A \in (2, \infty ),\; B \in [0,\,3/2]$
with
$\overline{\mathcal D}\subset B_{S_1}, \, $
$ A+\min\{B,1\}>3,\;  A+B\ge7/2 $.
Take $F: {\mathbb R}^3\mapsto\mathbb R^3 $ measurable with
$F|{B_{S_1}}\in L^{p_0}(B_{S_1})^3$,
$$|F(y)|\le\gamma\cdot|y|^{-A}\cdot s(y)^{-B}\ \hbox{for}\ y\in B^c_{S_1}.$$
\noindent
Let $u\in L^6(\overline{ \mathcal D}^c)^3\cap W^{1,1}_{loc}(\overline{ \mathcal D}^c)^3,
\pi\in L^2_{loc}(\overline{ \mathcal D}^c)$
with
$\nabla u \in L^2(\overline{ \mathcal D}^c)^9, \hbox{\rm div}\,u=0 $ and
\begin{eqnarray*}
\displaystyle{\int_{\overline{ \mathcal D}^c}
\left[\nabla u\cdot\nabla \varphi+((\tau e_{1}-\rho e_1\times z )\cdot\nabla u+\rho
e_1\times u
\qquad\qquad\qquad\qquad\right.} & \\
\displaystyle{ \left.+\tau(u\cdot\nabla)u-F)\cdot\varphi-\pi\,\hbox{\rm div}\,\varphi \right] \,
\hbox{\rm d}x=0
}
\end{eqnarray*}
for $\varphi\in C_0^{\infty}(\overline{ \mathcal D}^c)^3.$ Let $S\in(S_1,\infty).$  Then
\begin{eqnarray} \label{AA}
 \displaystyle{|\partial^\alpha u(x)|\le C\,
 (|x|s(x))}^{-1-|\alpha|/2} \ \ \hbox{for}\ x\in B^c _S,\ \alpha\in\mathbb N^3_0
\ \hbox{with}\ |\alpha|\le1
%,
.
\end{eqnarray}
\end{thm}

\bigskip
Now, using
Theorems  \ref{thm_Main_1} and \ref{nonlin},
we are in the position to prove the result on the decay of the pressure in the non-linear case:\

\begin{thm}\label{theorem2.7}
Consider the situation in Theorem \ref{nonlin}. Suppose in addition that
$A\ge 5/2$.
Then there is $c_0 \in \mathbb{R} $ such that inequality (\ref{CC}) holds.
\end{thm}
{\it Proof:}
Observe that $(u \cdot \nabla )u \in L^{3/2}( \overline{ \mbox{$\mathcal D$} }^c)^3$.
Thus, putting $p :=\min\{3/2,\, p_0\},\; \widetilde{ F}:=F-\tau \,(u \cdot \nabla )u,$
we get $\widetilde{ F}|\mathcal D_{S_1}\in L^p(\mathcal D_{S_1})^3$.
Put $B ^{\prime} := \min\{5/2,\, A+B-5/2\}$. Since $A\ge 5/2$, we have
\begin{eqnarray*}
|F(z)|
\le
\gamma \,C(S_1,A)\, |z|^{-5/2}\, s(z)^{-B ^{\prime} }
\quad \mbox{for}\;\;
z \in B_{S_1}^c.
\end{eqnarray*}
On the other hand, by Theorem \ref{nonlin} with $(S_1+S)/2$ in the place of $S$,
\begin{eqnarray*} &&
| \bigl(\, u(z)  \cdot \nabla )u(z) \,\bigr) |
\le
C\, |z|^{-5/2}\,s(z)^{-5/2}
\le
C\, |z|^{-5/2}\,s(z)^{-B ^{\prime} }
\end{eqnarray*}
for $ z \in B_{(S_1+S)/2}^c.$
In this way we get
$
| \widetilde{ F}|
\le
C\, |z|^{-5/2}\, s(z)^{-B ^{\prime} }
$
for
$
z \in B_{(S_1+S)/2}^c
$.

We further note that $B ^{\prime} \in (1/2,\, \infty )$. This is obvious
in the case $B ^{\prime} = 5/2$. If $B ^{\prime} < 5/2$, we have
$B ^{\prime} = A+B-5/2$. Due to the assumption $A+\min\{1,B\}> 3$ in Theorem \ref{nonlin},
we thus get $B ^{\prime} \in (1/2, \, \infty )$. (The requirement $A+B\ge 7/2$ in
Theorem \ref{nonlin} even yields $B ^{\prime} \ge 1$, but if this requirement is
weakened in a suitable way, pointwise decay of $u$ and $\nabla u$ could still be
proved. However, this point is not elaborated in \cite{DKN5}, and therefore is not reflected in
Theorem \ref{nonlin}. But we still take account of it here by avoiding to use the assumption
$A+B\ge 7/2$.)

We further have $u \in W^{1,p}_{loc}( \overline{ \mbox{$\mathcal D$} }^c)^3,\;
\pi \in
\color{black}
L^p_{loc}(\overline{ \mbox{$\mathcal D$} }^c)
\color{black}
,$ and
equation (\ref{2.**}) holds with $F$ replaced by $\widetilde{ F}$. Since in addition
$u|B_S^c \in L^6(B_S^c)^3,\; \nabla u|B_S^c \in L^2(B_S^c)^9$ and ${\rm div}\, u =0,$
we see that the assumptions of Theorem \ref{thm_Main_1} are satisfied with $p$ as defined
above and with $(S_1+S)/2,\; B ^{\prime} ,\, \widetilde{ F}$ in the role of $S_1,\,B$ and $F$, respectively.
Thus Theorem \ref{thm_Main_1} implies the conclusion of Theorem \ref{theorem2.7}.
\hfill $\Box $

\bigskip

\section{Formulation of the problem with artificial boundary conditions }

Recall that we defined $\mathcal D _R =  B_R\setminus\overline{\mathcal D}$.
We introduce the subspace $W_R$ of $H^1(\mathcal D_R)$ denoting
$$
W_R := \{v\in H^1(\mathcal D_R)^3 : v|{\partial \mathcal D} = 0\},
$$
where $v|\partial\mathcal D$ means the trace of $v$ on $\partial\mathcal D.$
\begin{lem}\label{lem_2.1}
{\rm
\color{black}
(\cite[Lemma 4.1]{DeFEM})
\color{black}
}
The estimate
$$
\|u\|_2 \le C \, (R\, \left\|\nabla u\right\|_2 + R^{1/2}\, \left\|u|{\partial B_R}\right\|_2)
$$
holds
for $R\in (0,\infty)$ with $\overline{\mathcal D} \subset B_R$ and for $u\in W_R$.
\end{lem}

\rm
We introduce an inner product $(\cdot,\,\cdot )^{(R)}$ in $W_R$ by defining
\begin{eqnarray*}
(v,w)^{(R)}= &\int_{\mathcal D_R}
\nabla v \cdot \nabla w
\, dx
%\\&&
+ \int_{\partial B_R} ( \tau /2)v\cdot w \,do_x \ \text{for} \ \  v,w \in W_R.
\end{eqnarray*}
The space $W_R$ equipped with this inner product is a Hilbert space.
The norm generated by this scalar product $(\cdot,\cdot)^{(R)}$ is denoted by $|\cdot|^{(R)}$, that is
$$
|v|^{(R)} : = \left(\left\|\nabla v\right\|^2_2+( \tau /2)\left\|v{|\partial B_R}\right\|^2_2\right)^{1/2} \ \text{for} \ v\in W_{R}.
$$
We define the bilinear forms
\begin{align*}
&a_R :H^1(\mathcal D_R)^3\times H^1(\mathcal D_R)^3\rightarrow \mathbb R,\\
&\beta_R : H^1(\mathcal D_R)^3 \times L^2(\mathcal D_R) \rightarrow \mathbb R,\\&\delta_R:
H^1(\mathcal D_R)^3\times H^1(\mathcal D_R)^3\rightarrow \mathbb R,\\
a_R (u,w):=&\int_{\mathcal D_R} [\nabla u\cdot \nabla w+ \tau  \partial_1 u\cdot w ] dx\\
&
+\frac{\tau }{2} \int_{\partial B_R} (u(x) \cdot w(x)) \left(1-\frac{x_1}{R}\right) do_x,\\
\beta_R (w,\sigma):= &- \int_{\mathcal D_R} ({\rm div} \, w) \, \sigma dx,
\\
\delta_R (u,w):=&
\color{black}
\int_{\mathcal D_R} \bigl[\, -\bigl(\, ( \varrho e_1\times x) \cdot \nabla \,\bigr)  u(x)
+\bigl(\,  \varrho e_1\times u(x) \,\bigr) \,\bigr] \cdot w(x)\, dx
\color{black}
\\ \hbox{for} \ \
u,w \in H^1(\mathcal D_R)^3,& \ \sigma \in L^2(\mathcal D_R),\
R\in (0,\infty)
\; \mbox{with}\; \overline{ \mathcal{D}}\subset B_R.
\end{align*}

\begin{lem}\label{lem_2.3}

{\it
Let $R \in (0, \infty ) $ with $\overline{ \mathcal{D}}\subset B_R$. Then
\begin{align*}
|a_R (u,w) + \delta _R(u,w)|
\le& C(R)\, |u|^{(R)}\,|w|^{(R)}
\end{align*}
for $ u,w \in H^1(\mathcal D_R)^3$.
}
\end{lem}
{\it Proof:} The proof of  Lemma \ref{lem_2.3} is based on use of Lemma \ref{lem_2.1}.

\medskip
The key observation in this section is stated in the following lemma, which is the basis of the theory
presented in this section.
\begin{lem}\label{lem_2.4}

{\it
Let $R\in (0, \infty ) $ with $\overline{ \mathcal{D}}\subset B_R$,
and let
$w\in W_R$. Then the equation $(|w|^{(R)})^2 = a_R(w,w)+\delta _R(w,w)$ holds.
}
\end{lem}
\smallskip

\noindent
{\it Proof:} Using the definitions $a_R(\cdot,\cdot),\, \delta_R(\cdot,\cdot)$, we get
\begin{eqnarray*} &&
a_R(w,w)+\delta_R(w,w)
\\&&
= \int_{\mathcal D_R} \left[|\nabla w|^2 + \tau
\partial _1 \left(\frac{|w|^2}{2}\right)
- ( \varrho e_1\times x) \cdot \nabla \left(\frac{|w|^2}{2}\right)\right] dx
\\&&
+\frac{\tau }{2} \int_{\partial B_R} |w(x)|^2 \left(1-\frac{x_1}{R}\right)do_x
\\&&
=\int_{\mathcal D_R} |\nabla w|^2\, dx
+ \int_{\partial B_R}\left(\frac{\tau }{2} |w(x)|^2
\frac{x_1}{R} - \frac12 ( \varrho e_1\times x)\cdot \frac{x}{R} |w(x)|^2\right)do_x
\\&&
+\frac{\tau }{2} \int_{\partial B_R} |w(x)|^2\left(1-\frac{x_1}{R} \right) do_x
\\&&
=\int_{\mathcal D_R}|\nabla w|^2\, dx + \frac{\tau }{2} \int_{\partial B_R} |w(x)|^2 = (|w|^{(R)})^2.
\end{eqnarray*}
We applied that
\color{black}
$$( \omega \times x)\cdot x =0 \; \mbox{for}\;\; x,\, \omega \in \mathbb{R}^3 .$$
\color{black}\

As in \cite{DeKr},
we obtain that the bilinear form $\beta_R$ is stable:

\begin{thm}\label{thm_2.2}
{\rm (\cite[Corollary 4.3]{DeKr})}
{\it
Let $R>0$ with $\overline{ \mathcal D}\subset B_R$.
Then
\[
\inf\limits_{\rho \in L^2(\mathcal D_R), \rho \not=0}
\sup_{v\in W_R, v\not=0} \frac{\beta_R(v,\rho)}{|v|^{(R)}\|\rho\|_2} \ge C(R).
\]
}
\end{thm}
\color{black}
We note that functions from $W^{1,1}_{loc}( \overline{ \mathcal D}^c)$
with $L^2$-integrable gradient are
$L^2$-integrable on truncated exterior domains:
\begin{lem}[\mbox{\cite[Lemma II.6.1]{Galdi}}] \label{lemma3.xx}
Let $ w \in W^{1,1}_{loc}( \overline{ \mathcal D}^c)$
with $\nabla w \in L^2( \overline{ \mathcal D}^c)^3$, and let $R \in (0, \infty ) $
with $\overline{ \mathcal D}\subset B_R.$ Then
$w{| \mathcal D_R}
\in L^2(\mathcal D_R).
$
In particular the trace of $w$ on $\partial \mathcal D$ is well defined.
\end{lem}
The preceding lemma is implicitly used in the ensuing theorem, where
we introduce
an
extension operator
$\mbox{$\mathfrak E$}: H ^{1/2} ( \partial \mathcal{D})^3 \mapsto W^{1,1}_{loc}( \overline{ \mathcal{D}}^c)^3$
such that
$\mbox{div}\, \mbox{$\mathfrak E$} (b)=0
$.
\begin{thm}[\mbox{\cite[Exercise III.3.8]{Galdi}}]\label{thm_2.1}
{\it
There is an operator
$\mbox{$\mathfrak E$}$
from
$H^{1/2}(\partial \mathcal D)^3 $
into
$W^{1,1}_{loc}( \overline{ \mathcal{D}}^c)^3$
satisfying the relations
$
\nabla \mbox{$\mathfrak E$} (b) \in L^2( \overline{ \mathcal{D}}^c)^9,\;
\mbox{$\mathfrak E$} (b)|{\partial \mathcal{D}} = b$
and
${\rm div}\,\mbox{$\mathfrak E$} (b) =0
$
for
$
b\in H^{1/2}(\partial \mathcal{D})^3.
$
}
\end{thm}
\color{black}

\vspace{1ex}
In view of Lemma \ref{lem_2.3} and \ref{lem_2.4}
and Theorem \ref{thm_2.1} and \ref{thm_2.2}, the theory of mixed variational problems
yields
\begin{thm}\label{thm_2.3}
{\it
Let $S>0$ with $\overline{ \mathcal D}\subset B_S,\;
R \in [2S, \infty ),\;
F \in
L^{6/5}( \mathcal{D}_R)^3, \;
b \in H ^{1/2} ( \partial \mathcal{D})^3.$
Then there is a uniquely determined pair of functions
$(\widetilde V, P )= \bigl(\, \widetilde V(R,F,b),\, P(R,F,b) \,\bigr) \in W_R \times L^2( \mathcal{D}_R)$
such that
\begin{eqnarray} &&
a_R( \widetilde V,g) + \delta _R( \widetilde V,g) +\beta _R(g, P )
\label{eq_2_1}\\&&
=
\int_{ \mathcal{D}_R}F \cdot g\, dx - a_R \bigl(\, \mbox{$\mathfrak E$} (b)|\mathcal{D}_R,\, g \,\bigr)
-\delta _R \bigl(\, \mbox{$\mathfrak E$} (b)|\mathcal{D}_R,\, g \,\bigr)
\;\; \mbox{for}\;\; g \in W_R,
\nonumber\\[1ex]&&
\beta _R( \widetilde V, \sigma )=0
\;\; \mbox{for}\;\; \sigma \in L^2( \mathcal{D}_R),
\label{eq_2_2}\end{eqnarray}
where the operator \mbox{$\mathfrak E$} was introduced in Theorem \ref{thm_2.1}.
}
\end{thm}

Let us interpret variational problem (\ref{eq_2_1}), (\ref{eq_2_2})
as a boundary value problem.
Define the expression used in the boundary condition on the artificial boundary $\partial B_R:$\
$$\mathcal L _R (u,\pi)(x):=  \left(\sum^3_{j=1} \partial_j u_k(x) \frac{x_j}{R}
\color{black}
- \pi(x) \frac{x_k}{R}
\color{black}
+ \frac{\tau }{2}
\left(1-\frac{x_1}{R}\right) u_k(x)\right)_{1\le k\le 3}$$
for $x \in \partial B_R,\; R \in (0, \infty ) $ with $\overline{ \mbox{$\mathcal D$} } \subset B_R,\;
\color{black}
u\in W^{2,\, 6/5}(\mathcal D_R)^3,
\color{black}
\ \pi\in W^{1,\, 6/5}(\mathcal D_R)$.

\begin{lem}\label{lem_interpret}  Assume that $ \mathcal D$ is $\mathcal C^2$-bounded. Let
$
\color{black}
S \in (0, \infty )
\color{black}
$  with
$\overline{ \mathcal D}\subset B_S,\;
\color{black}
R\in[2S, \infty),
\color{black}
\;
F \in L^{6/5}( \mathcal{D}_R)^3$ and
$b\in W^{7/6,\,6/5 }( \partial\mathcal{D})^3.$
Put $V:=\widetilde V(R,F,b)+\mathfrak E(b)|\mathcal D_R$,
with $V(R,F,b)$ from Theorem \ref{thm_2.3} and $\mbox{$\mathfrak E$} (b)$
from Theorem \ref{thm_2.1}.
Suppose that
$
V\in W^{2,6/5}(\mathcal D_R)^3$
and
$ P=P(R,F,b)\in
W^{1,\, 6/5}(\mathcal D_R)$, with $P(R,F,b)$ also introduced in Theorem \ref{thm_2.3}.
Then
\begin{equation}
\begin{array}{crl}
-\Delta V (z) + ( \tau e_1- \varrho e_1 \times z)\cdot \nabla V(z) + \varrho e_1 \times V (z)+\nabla P(z)
= F(z),
\\
{\rm div}\,V(z)=0 \
\end{array}
\end{equation}

\noindent
for $z\in \mathcal D_R,$ and
$V|\partial\mathcal D =b, \ \ \mathcal L_R(V,P)=0.$
\end{lem}

The proof of  Lemma \ref{lem_interpret} is obvious. This lemma
means that a solution of variational problem  (\ref{eq_2_1}), (\ref{eq_2_2}) may be considered as a weak solution of the modified Oseen system with rotation in $\mathcal D_R$, under the Dirichlet boundary condition on $\partial \mathcal D$ and under the artificial boundary condition $\mathcal L_{R}(V,P)=0 $ on $\partial B_R$. The solution of   (\ref{eq_2_1}), (\ref{eq_2_2}) will be now compared to the exterior modified Oseen flow
introduced in
Corollary \ref{corollary2.1}:

\begin{thm}\label{thm_Main_2}
Suppose that \mbox{$\mathcal D$} is $C^2$-bounded. Let $\gamma ,\, S_1 \in (0, \infty ) $
with $\overline{ \mathcal D} \subset B_{S_1}, \; A \in [5/2,\, \infty ),\; B \in \mathbb{R} $
with
$\ A+\min\{1,B\} >  3$.
Let $F:\overline{\mathcal D}^c \mapsto \mathbb{R}^3 $ be measurable with
$F| \mathcal D_{S_1} \in L^{6/5}(\mathcal D_{S_1})^3$ and
$|F(z)|\le \gamma\,|z|^{-A}s(z)^{-B} \ \text{for} \ z\in B^c_{S_1}$.

Let $ b \in W^{7/6,\, 6/5}( \partial \mbox{$\mathcal D$} )^3,\;
u \in W^{1,1}_{loc}( \overline{ \mathcal D}^c)^3\cap L^6(\overline{ \mathcal D}^c)^3$
such that
$
\nabla u \in L^2(\overline{ \mathcal D}^c)^9,\; {\rm div}\, u =0,\;
u| \partial \mbox{$\mathcal D$} =b
$
and equation (\ref{2.*}) is satisfied.

For $R \in [2S_1,\, \infty ),$ put
$
V_R:=\widetilde V(R,F,b) + \mbox{$\mathfrak E$} (b),%\; P := P(R,F,b),
$
with $\mbox{$\mathfrak E$} (b)$ from Theorem \ref{thm_2.1},
and
$\widetilde V(R,F,b)
$
from Theorem \ref{thm_2.3}.
Then
\begin{align*}
&|u|_{\mathcal D_R}-V_R|^{(R)} \le C\,R^{-1} \quad \mbox{for}\;\; R \in [2S, \infty ).
\end{align*}
\end{thm}
We note that since $W^{2,\, 6/5}( \mbox{$\mathcal D$} ) \subset H^{1}( \mbox{$\mathcal D$} )$
by a Sobolev inequality, we have
$W^{7/6,\, 6/5}( \partial \mbox{$\mathcal D$} ) \subset H^{1/2}( \partial \mbox{$\mathcal D$} )$,
as follows with the usual
\color{black}
lifting and trace properties.
\color{black}
As a consequence,
$b \in H ^{1/2} ( \partial \mbox{$\mathcal D$})^3$, so the term $\mbox{$\mathfrak E$} (b)$
is well defined. We further remark that by Corollary \ref{corollary2.1} with $p=6/5$,
the function $F$ may be considered as a bounded linear functional on
$\mbox{$\mathcal D$} ^{1,2}_0( \overline{ \mbox{$\mathcal D$} }^c)^3.$
Therefore, as explained in Remark \ref{rem_existence}, a function $u$ with properties as stated in
Theorem \ref{thm_Main_2} does in fact exist.

\smallskip
\noindent
{\it Proof of Theorem \ref{thm_Main_2}:}
All conditions in Corollary \ref{corollary2.1} are verified if $\gamma ,\, S_1,\, A,\, B,\, F,\, u$
are given as in Theorem \ref{thm_Main_2}, and if $p=6/5$ and $S=2\,S_1.$
Note in this respect that the conditions on $u$ in Theorem \ref{thm_Main_2} obviously imply
$u \in W^{1,\, 6/5}_{loc}( \overline{ \mbox{$\mathcal D$} }^c)^3$.
Corollary \ref{corollary2.1} now yields that
$F \in L^{6/5}(\overline{ \mbox{$\mathcal D$} }^c)^3$ and that the function $u$
satisfies inequalities (\ref{1.6}) and (\ref{1.7})
with $S=2\, S_1$.

On the other hand, since
$u \in W^{1,\, 6/5}_{loc}( \overline{ \mbox{$\mathcal D$} }^c)^3$,
the function $G$ already considered in the proof of Corollary \ref{corollary2.1} (see (\ref{BB}))
belongs to $L^{6/5}_{loc}( \overline{ \mbox{$\mathcal D$} }^c)^3$.
Therefore, by interior regularity of solutions to the Stokes system (see \cite[Theorem IV.4.1]{Galdi}),
we may deduce from the equations (\ref{2.*}) and ${\rm div}\, u =0$ that
$u \in W^{2,\, 6/5}_{loc}( \overline{ \mbox{$\mathcal D$} }^c)^3$
and that there is $\pi \in W^{1,\, 6/5}_{loc}( \overline{ \mbox{$\mathcal D$} }^c)^3$ with
$L(u)+\nabla \pi =F$. In particular the pair $(u,\pi )$ verifies (\ref{2.**}).
In view of our assumptions on $A$ and $B$, we thus see that the requirements in
Corollary \ref{corollary2.2}
are fulfilled
for $\gamma ,\, S_1,\, A,\, B,\, F,\, u$
as in Theorem \ref{thm_Main_2} and for $p=6/5$ and $S=2\,S_1$.
As a consequence, Corollary \ref{corollary2.2} yields that there is $c_0 \in \mathbb{R} $ such that
(\ref{CC}) holds with $S=2\, S_1$.

Take $R \in [2\, S_1, \, \infty ).$
Since
$u \in W^{2,6/5}_{loc}( \overline{ \mbox{$\mathcal D$} }^c)^3$,
we have
$
\color{black}
u | \partial B_R \in W^{7/5,\,6/5}( \partial B_R)^3
\color{black}
$.
Combining this relation with the assumption
$
\color{black}
b
\in W^{7/5,\,6/5}( \partial \mbox{$\mathcal D$})^3
\color{black}
$
and the boundary condition $u| \partial \mbox{$\mathcal D$} =b,$
we get $u| \partial \mbox{$\mathcal D$} _R \in W^{7/5,\,6/5}( \partial \mbox{$\mathcal D$}_R)^3$.
Moreover our requirements on $u$ yield that
$u| \mbox{$\mathcal D$} _R \in W^{1,\,6/5}( \mbox{$\mathcal D$}_R)^3$.
Since $F \in L^{6/5}( \overline{ \mbox{$\mathcal D$} }^c)^3$,
as already mentioned, we get
$G| \mbox{$\mathcal D$} _R \in L^{6/5}( \mbox{$\mathcal D$}_R)^3$,
with $G$ from (\ref{BB}).
Recalling that \mbox{$\mathcal D$} is supposed to be $C^2$-bounded,
we may now apply the result in \cite[Lemma IV.6.1]{Galdi} on boundary regularity
of solutions to the Stokes system. This reference yields that
$u| \mbox{$\mathcal D$} _R \in W^{2,\,6/5}( \mbox{$\mathcal D$}_R)^3,\;
\pi | \mbox{$\mathcal D$} _R \in W^{1,\,6/5}( \mbox{$\mathcal D$}_R)$
and that the pair $(u,\pi) $ solves (\ref{eq_1.1}).

\medskip
Let $P_R:=P(R,F,b)$ be given as in Theorem \ref{thm_2.3},
and put
$
w:=u-V_R,\; \kappa := \pi - P_R\, $,
and let  $\, g \in W_R$.
Note that by Theorem \ref{thm_2.3}, we have

$$
a_R(V_R,g)+\delta _R(V_R,g) + \beta _R(g,P_R) = \int_{ \mathcal D_R}F \cdot g\, dx.
$$

Thus
\begin{eqnarray*} &&\hspace{-2em}
a_R (w,g) +\delta_R(w,g) + \beta _R (g,\kappa)
\\&&\hspace{-2em}
=a_{R} (u|_{\mathcal D_R},g) +\delta _R(u|_{\mathcal D_R},g)
 + \beta_R(g,\pi|_{\mathcal D_R})- \bigl(\, \underbrace{a_R (V_R ,g )
+  \delta _R (V_R,g) +\beta_R(g, P_R )\bigr)}
\\&&\hspace{25em}
_{=\int_{\mathcal D_R} F\cdot g\, dx} \,
\end{eqnarray*}
\begin{eqnarray*}
&&\hspace{-2em}=\int_{\mathcal D_R} \bigl(\, \nabla u \cdot \nabla g + \tau  \partial_1 u \cdot g
- ( \varrho e_1 \times x) \cdot \nabla u\cdot g
+( \varrho e_1 \times u) \cdot g
\\&&
\hspace{18em}- \pi\,{\rm div} \, g  - F \cdot g \,\bigr) \, dx
\\&&\hspace{-1em} +\frac{\tau }{2} \int_{\partial B_R} u(x) \cdot  g(x) \left(1-\frac{x_1}{R}\right) do_x
%\end{eqnarray*}
%
%\begin{eqnarray*}
\\&&\hspace{-2em}
=\displaystyle\int_{\mathcal D_R} [
%\underbrace{
-\Delta u \cdot g + \tau  \partial_1 u \cdot g
- ( \varrho e_1\times x) \cdot  \nabla u\cdot g
+( \varrho e_1 \times u)  \cdot  g
\\&&\hspace{18em}+ \nabla  \pi \cdot g - F\cdot g
%}_{=0}
]\, dx
\\&&
\underbrace{\hspace{-1EM}+\int_{\partial B_R}\Bigl(
 \sum^3_{j,k=1} \bigl[\, \partial_j u_k(x)\, g_k(x) \frac{x_j}{R}
-\pi (x)\delta_{jk} g_k(x) \frac{x_j}{R}  \,\bigr]+ \frac{\tau }{2} u(x) \cdot g(x) (1-\frac{x_1}{R})\Bigr)  do_x.}
\\&&
\hspace{15em}
_{=\int_{\partial B_R}\mathcal L_R(u,\pi)\cdot g \ do}
\end{eqnarray*}

\noindent

\medskip
\noindent
Since the pair $(u,\pi)$ solves (\ref{eq_1.1}),
we  now get
\begin{eqnarray}
a_R (w,g) +\delta_R(w,g) + \beta _R (g,\kappa)
=
\int_{\partial B_R}
\mathcal L_R  (u,\pi)(x) \cdot  g(x)\, do_x.
\end{eqnarray}
Let $c\in\mathbb R$ be an arbitrary constant.
For $g:=w$ we get with Lemma \ref{lem_2.4} that
\begin{eqnarray}\label{eq_w}
(|w|^{(R)})^2 &=&a_R (w,w) +\delta_R(w,w) + \beta _R (w,\kappa)\nonumber\\ &=&
\int_{\partial B_R}
\mathcal L_R  (u,\pi+c)(x) \cdot  w(x)\, do_x\, ,
\end{eqnarray}
\color{black}
because by the assumptions on $u$ and Theorem \ref{thm_2.1} and \ref{thm_2.3},
\color{black}
$$
\int_{\partial B_R}\Bigl[
 \sum^3_{j,k=1} \,  c\delta_{jk} w_k(x) \frac{x_j}{R}  \,\Bigr]  do_x=
 \int_{
\color{black}
\partial {\mathcal D}
\color{black}
}
  \,  c\,w\cdot n   \,  do_x+\int_{ {\mathcal D}_R}
  \,  c\,\hbox{div}\, w \,  d x=0,
$$
\color{black}
where $n$ denotes the outward unit normal to
$\mathcal D.$
\color{black}
Let $c_0$ be the constant introduced above as part of
estimate (\ref{CC}).
Because
 $$
\int_{\partial B_R} \mathcal L_R  (u,\pi+c_0)(x) \cdot  w(x)\, do_x
\le \|\mathcal L_R(u,\pi+c_{0})\|_2\,
\|w|_{
\color{black}
\partial B_R
\color{black}
} \|_2
\le C\, \|\mathcal L_R (u,\pi+c_0)\|_2 \cdot |w|^{(R)},
$$
we get from (\ref{eq_w}) $$|w|^{(R)}\le C \|\mathcal L_R (u,\pi+c_0)\|_2.
$$
The last step is estimation:
 $\|\mathcal L_R (u,\pi+c_{0}) \|_2\le
C\,\cdot
R^{-1}$.
\color{black}
We start by observing that
\color{black}
\begin{eqnarray*}&&
\hspace{-4em}\|\mathcal L_R (u,\pi+c_{0}) \|_2
\\&&
\hspace{-2em}\le C \Big[\|\nabla u|_{\partial B_R} \|_2 + \|[\pi(x)+c_0] |_{\partial B_R}\|_2 +
\left(\int_{\partial B_R} \left(1-\frac{x_1}{R}\right)^2|u(x)|^2do_x\right)^{1/2}\Big].
\end{eqnarray*}
As explained above, inequalities (\ref{1.6}), (\ref{1.7}) and (\ref{CC})
are valid with $S=2\, S_1$. According to (\ref{1.6}) and (\ref{CC}), we have
$
| u(x)| \le C\, (|x|s(x))^{-1},
$
and $|\pi(x)+c_0|\le C\,|x|^{-2} $
for $x \in B_{2\cdot S_1}^c$.
Inequality (\ref{1.7}) yields
$
| \nabla  u(x)| \le C\,
|x|^{
\color{black}
-3/2
\color{black}
}
s(x)^{-B ^{\prime} }
$
for $x$ as before, with
$B ^{\prime} :=3/2-\max\{0,\, 7/2-A-B\}$.
If $B\ge 1$, we recall that $A\ge 5/2$, getting $A+B\ge 7/2,$
hence $B ^{\prime} =3/2$.
On the other hand, if $B<1$, then
$
\color{black}
\min
\color{black}
\{1,B\}=B$, so that the assumption
$A+
\color{black}
\min
\color{black}
\{1,B\}>3$
becomes $A+B>3$, hence $B ^{\prime} > 1$. Thus we get in any case that $B ^{\prime} >1>1/2$.
In view of these observations, and with Lemma \ref{Fa}, we obtain
\begin{eqnarray*}
\|\mathcal L_R(u,\pi +c_0)\|_2 &\le& C
\left[ \left(
\int_{\partial B_R}\, |x|^{-3}s(x)^{-2\,B ^{\prime} } \,  do_x
\right)^{1/2}
+ \left(\int_{\partial B_R} |\pi(x)+c_0|^2 do_x\right)^{1/2}\right.
\\&&\hspace{3em}
+\left. \left(\frac1{R^2} \int_{\partial B_R}(|x|-x_1)^2|u(x)|^2do_x\right)^{1/2}\right]
\\
&\le& C \left[ \left(\frac1{R^3}\int_{\partial B_R}s(x)^{-2\, B ^{\prime} } do_x\right)^{1/2}
+ \left(\frac{1}{R^4}\int_{\partial B_R} 1\, do_x\right)^{1/2}
\right.
\\&&\hspace{3em}
\left.+\left(\frac1{R^2}\int_{\partial B_R} s(x)^2(|x|s(x))^{-2}do_x\right)^{1/2}      \right]
\\
&\le&
C \left[\left(\frac{1}{R^2} \right)^{1/2}
+\left(\frac{1}{R^2} \right)^{1/2}  +\left(\frac1{R^4}\int_{\partial B_R}1\, do_x\right)^{1/2}\right]\le CR^{-1}.
\end{eqnarray*}
This completes the proof of Theorem \ref{thm_Main_2}. \hfill $\Box $

\bigskip

\bigskip

\textbf{Acknowledgements:} \vskip0.25cm \textit{ The works of S.K. and \v{S}. N. were supported by Grant No. 16-03230S of GA\v{C}R in the framework of RVO 67985840, S.K. is supported by RVO 12000. Final version was supported by Grant No. 19-04243S of GA\v CR.}

\end{document}